\documentclass[12pt,normalheadings]{scrartcl}

\usepackage[utf8]{inputenc}
\usepackage[ngerman]{babel}
\usepackage{amsmath,amssymb}
\usepackage[thmmarks]{ntheorem}
\usepackage{verbatim}
\usepackage{graphicx}
\usepackage[all]{xy}

\newcommand{\otr}{\underset{R}{\otimes}}
\newcommand{\Ass}{\operatorname{Ass}}
\newcommand{\dotcup}{\mathrel{\mathaccent\cdot\cup}}

\theorembodyfont{\itshape}
\theoremstyle{plain}
\newtheorem{Satz}{Satz}[section]
\newtheorem{Lemma}[Satz]{Lemma}
\newtheorem{Folgerung}[Satz]{Folgerung}
\theorembodyfont{\upshape}
\newtheorem{Bemerkung}[Satz]{Bemerkung}

\theoremheaderfont{\sc}
\theorembodyfont{\upshape}
\theoremstyle{nonumberplain}
\theoremseparator{}
\theoremsymbol{$\Box$}
\newtheorem{Beweis}{Beweis}

\title{\Large Totalseparierte Moduln}
\author{\large Helmut Zöschinger\\
  \large Mathematisches Institut der Universität München\\
  \large Theresienstr. 39, D-80333 München\\
  \large E-mail: zoeschinger$@$mathematik.uni-muenchen.de
}
\date{}

\newcounter{abccount}
\newenvironment{abc}{%
\begin{list}{\alph{abccount})}{%
\usecounter{abccount}%
\setlength{\partopsep}{0pt}%
\setlength{\topsep}{1ex}%
\setlength{\itemsep}{0pt}%
}%
}{\end{list}}

\newcounter{iiicount}
\newenvironment{iii}{%
\begin{list}{\roman{iiicount})}{%
\usecounter{iiicount}%
\setlength{\labelwidth}{3em}%
\setlength{\partopsep}{0pt}%
\setlength{\topsep}{1ex}%
\setlength{\itemsep}{0pt}%
}%
}{\end{list}}

\begin{document}
\maketitle

\centerline{\textbf{Abstract}}
\begin{abstract}
  \noindent
  Let $(R, \mathfrak{m})$ be a noetherian local ring, $M$ a separated
  $R$-module (i.e. $\bigcap\limits_{n\geq 1}\mathfrak{m}^n M = 0$) and
  $\widehat{M} = \lim\limits_{\leftarrow} M/\mathfrak{m}^n M$ its
  completion. Generally, $M$ is not pure in $\widehat{M}$ and $\widehat{M}$
  is not pure-injective. But if $M$ is \emph{totally separated}, i.e. $X
  \otr M$ is separated for all finitely generated $R$-modules $X$, the
  situation improves: In this case, $M$ is pure in $\widehat{M}$ and, under
  additional conditions, $\widehat{M}$ is even pure-injective, e.g. if $M
  \cong X^{(I)}$ holds with $X$ finitely generated or $M \cong
  \coprod_{i=1}^{\infty} R/\mathfrak{m}^i$.\\ In section~2, we investigate
  the question under which conditions both $M$ and $\widehat{M}$ are totally
  separated and establish a close connection to the class of strictly
  pure-essential extensions. In section 3, we replace the completion
  $\widehat{M}$ in the case $M = \coprod_{i\in I}M_i$ with the
  $\mathfrak{m}$-adic closure $A$ of $M$ in $P = \prod_{i\in I} M_i$, i.e.
  with $A = \bigcap_{n \geq 1}(M + \mathfrak{m}^n P)$. We give criteria so
  that $A/M$ is radical and show that this always holds in the countable case
  $M = \coprod_{i=1}^{\infty} M_i$. Finally, we deal with the case that $A$
  is even totally separated and additionally determine the coassociated
  prime ideals of $A/M$.
\end{abstract}

\bigskip


\bigskip


\section{Wann ist $\widehat{M}$ rein-injektiv?}
Stets sei $(R, \mathfrak{m})$ ein noetherscher lokaler Ring. Ein $R$-Modul
$M$ heißt \emph{totalsepariert}, wenn für jeden endlich erzeugten $R$-Modul
$X$ auch $X \otr M$ separiert ist. Das ist nach (\cite{011} Lemma 2.1)
äquivalent damit, dass die rein-injektive Hülle $N$ von $M$ separiert ist,
oder auch damit, dass die kanonische Abbildung $M \to \prod_{n \geq 1} M /
\mathfrak{m}^n M$ ein reiner Monomorphismus ist. Die Klasse aller
totalseparierten $R$-Moduln ist gegenüber beliebigen direkten Produkten und
reinen Untermoduln abgeschlossen, und weil $\widehat{M}$ ein Modul zwischen
$M$ und $\prod_{n \geq 1} M / \mathfrak{m}^n M$ ist, ist jeder
totalseparierte $R$-Modul $M$ rein in $\widehat{M}$.

\begin{Lemma}\label{1.1}
Seien $M$ und $\widehat{M}$ totalsepariert und $X$ endlich erzeugt. Dann
induziert die Einbettung $M \subset \widehat{M}$ einen Isomorphismus
\begin{equation*}
  \widehat{X \otr M} \longrightarrow X \otr \widehat{M}\ .
\end{equation*}
\end{Lemma}

\begin{Beweis}
  Vorbemerkung: Ist $B$ ein separierter $R$-Modul und $A \subset B$ ein
  dichter Untermodul (d.\,h. $B/A$ radikalvoll), so ist $\widehat{A} \to
  \widehat{B}$ nach Simon (\cite{005} Lemma~1.2) surjektiv. War $A$
  zusätzlich ein Unterraum von $B$, ist also $\widehat{A} \to
  \widehat{B}$ ein Isomorphismus.

  In unserer Situation ist $M$ rein und dicht in $\widehat{M}$, also auch $X
  \otr M$ rein und dicht in $X \otr \widehat{M}$, und weil nach
  Voraussetzung $X \otr \widehat{M}$ separiert, also als Faktormodul von
  $R^n \otr \widehat{M}$ sogar vollständig ist, liefert die Vorbemerkung
  die Behauptung.
\end{Beweis}

\begin{Bemerkung}\label{1.2}
  Ist $M$ frei, so findet sich die Formel (\ref{1.1}) für zyklische $X$ bei
  Simon (\cite{006} p.~974), für endlich erzeugte $X$ bei Marley-Webb
  (\cite{003} p.~10).
\end{Bemerkung}

\begin{Folgerung}\label{1.3}
  Sei $M$ totalsepariert, $\widehat{M}$ rein-injektiv und $X$ endlich
  erzeugt. Dann ist auch $\widehat{X \otr M}$ rein-injektiv.
\end{Folgerung}

\begin{Beweis}
  Weil $\widehat{M}$ rein-injektiv und separiert ist, gibt es einen
  zerfallenden Monomorphismus $\widehat{M} \longrightarrow \prod_{\lambda
    \in \Lambda} A_{\lambda}$, in dem alle $A_{\lambda}$ von endlicher Länge
  sind (siehe \cite{011} Lemma~1.5,c). Dann ist auch $X \otr \left(\prod
    A_{\lambda}\right)$ rein-injektiv, also auch $X \otr \widehat{M}$ und
  mit ($\ref{1.1}$) folgt die Behauptung.
\end{Beweis}

Der Spezialfall $M = R^{(I)}$ interessiert uns besonders:

\begin{Satz}\label{1.4}
  Ist $X$ endlich erzeugt, so gilt für jede Indexmenge $I$, dass
  $\widehat{X^{(I)}}$ rein-injektiv ist.
\end{Satz}

\begin{Bemerkung}\label{1.5}
  Sind in $M = \coprod M_i$ alle $M_i$ endlich erzeugt (insbesondere alle
  $\widehat{M_i}$ rein-injektiv), so ist $\widehat{M}$ i.\,Allg. nicht
  rein-injektiv. In dem Beispiel von Jensen und Zimmermann-Huisgen
  (\cite{002} Lemma 3) ist $\mathfrak{m}^2 = 0$, $\dim_k
  (\mathfrak{m}/\mathfrak{m}^2) = 2$ und $\mathcal{R}$ ein
  Repräsentantensystem des Restklassenkörpers $k$. Mit $\mathfrak{m} =
  (x,y)$ und einer unendlichen Teilmenge $\Lambda \subset \mathcal{R}$ sind
  dann die Ideale $\mathfrak{a}_{\lambda} = (x-\lambda y)$ paarweise
  verschieden, die $M_{\lambda} = R/\mathfrak{a}_{\lambda}$ einreihig von
  der Länge 2, aber $M = \coprod_{\lambda\in\Lambda} M_{\lambda}$
  \emph{nicht} rein-injektiv, denn die rein-injektive Hülle von $M$ ist
  $\prod_{\lambda\in\Lambda} M_{\lambda}$.
\end{Bemerkung}

Damit also der kanonische Monomorphismus $\widehat{\coprod M_i}
\longrightarrow \widehat{\prod M_i} \cong \prod \widehat{M_i}$ zerfällt,
müssen gewisse Bedingungen an den Zusammenhalt der $M_i$ erfüllt sein: In
($\ref{1.4}$), dass alle $M_i \cong X$ sind mit $X$ endlich erzeugt; in
(\ref{1.10}), dass bei festem $n \geq 1$ fast alle $M_i / \mathfrak{m}^n
M_i$ als $R/\mathfrak{m}^n$-Moduln frei sind.

\begin{Lemma}\label{1.6}
  Sei $B$ ein separierter $R$-Modul und $A \subset B$ ein Untermodul, so dass
  für alle $n \geq 1$ gilt: Es ist $A \cap \mathfrak{m}^n B = \mathfrak{m}^n
  A$ und $B/A + \mathfrak{m}^nB$ als $R/\mathfrak{m}^n$-Modul frei. Dann ist
  $\widehat{A} \to \widehat{B}$ ein zerfallender Monomorphismus.
\end{Lemma}

\begin{Beweis}
  Wir zeigen im \textbf{1. Schritt} durch Induktion über $n \geq 1$, dass es
  Untermoduln $B \supset X_1 \supset X_2 \supset \ldots$ gibt mit $X_n
  \supset \mathfrak{m}^n B$ und
  \begin{equation*}
    \frac{X_n}{\mathfrak{m}^{n}B} \oplus \frac{A +
      \mathfrak{m}^{n}B}{\mathfrak{m}^{n}B} = \frac{B}{\mathfrak{m}^{n}B}
  \end{equation*}
  für alle $n \geq 1$. Klar ist $n = 1$, und hat man bei $n > 1$ bereits die
  $X_1, \dotsc, X_{n-1}$ wie gewünscht, folgt aus $X_{n-1} + A = B$, dass
  $X_{n-1}/ \mathfrak{m}^{n}B + (A + \mathfrak{m}^{n}B)/\mathfrak{m}^{n}B =
  B/\mathfrak{m}^{n}B$ ist, also mit Hilfe der zweiten Voraussetzung ein
  $X_{n-1} \supset X_n \supset \mathfrak{m}^{n}B$ existiert mit $X_n
  /\mathfrak{m}^{n}B \oplus (A+\mathfrak{m}^{n}B)/\mathfrak{m}^{n}B =
  B/\mathfrak{m}^{n}B$ wie verlangt.

  \bigskip

  Mit Hilfe dieser $X_n$ sei im \textbf{2. Schritt} $\rho_n \colon
  B/\mathfrak{m}^{n}B \to (A+\mathfrak{m}^{n}B)/\mathfrak{m}^{n}B$ die
  Projektion auf den zweiten Summanden, und weil $\omega_n \colon A /
  \mathfrak{m}^{n}A \to (A + \mathfrak{m}^{n}B)/\mathfrak{m}^{n}B$ nach der
  ersten Voraussetzung ein Isomorphismus ist, wird $\beta_n := \omega_n^{-1}
  \circ \rho_n$ ein Linksinverses von $\alpha_n \colon A/\mathfrak{m}^{n}A
  \to B/\mathfrak{m}^{n}B$. Wegen $X_n \supset X_{n+1}$ ist auch
  \begin{equation*}
    \xymatrix{
      A/\mathfrak{m}^{n}A & B/\mathfrak{m}^{n}B \ar[l]_{\beta_n} \\
      A/\mathfrak{m}^{n+1}A \ar[u] & B/\mathfrak{m}^{n+1}B
      \ar[l]_{\beta_{n+1}} \ar[u]
    }
  \end{equation*}
  kommutativ, so dass für den Homomorphismus
  $g=(\beta_1,\beta_2,\beta_3,\dotsc) \colon \prod_{n\geq 1}
  B/\mathfrak{m}^{n}B \to \prod_{n \geq 1} A/\mathfrak{m}^{n}A$ gilt
  $g(\widehat{B}) \subset \widehat{A}$ und der induzierte Homomorphismus
  $\widehat{g} \colon \widehat{B} \to \widehat{A}$ ein Linksinverses von
  $\widehat{A} \to \widehat{B}$ ist.
\end{Beweis}

\begin{Bemerkung}\label{1.7}
  Im Spezialfall $A = R^{(I)}$, $B = R^I$ ist die Aussage wohlbekannt
  (\cite{001} Proposition~3) und aus ihr folgt allgemeiner: Ist $M$ flach
  und separiert, so ist $\widehat{M}$ flach und rein-injektiv (\cite{001}
  Proposition~4).

  Auch wenn $B$ separiert und $B/A$ flach ist, sind die Voraussetzungen in
  (\ref{1.6}) erfüllt, denn $A$ ist dann rein in $B$ und $R/\mathfrak{m}^{n}
  \otr B/A \cong B/A + \mathfrak{m}^{n}B$ als $R/\mathfrak{m}^{n}$-Modul
  flach, also frei. Beispiel~1: Sind in $A = \coprod M_i$ alle $M_i$
  separiert und flach und ist $B =  \prod M_i$, so folgt $\widehat{A}
  \subset^{\oplus} \widehat{B}$. Beispiel~2 (\cite{007} p.~388): Ist $B$
  vollständig und $B/A$ separiert und flach, so folgt $A \subset^{\oplus} B$.
\end{Bemerkung}

\begin{Folgerung}\label{1.8}
  Seien in $M = \coprod M_i$ alle $M_i$ totalsepariert und flach und sei $X$
  endlich erzeugt. Mit $P = \prod M_i$ ist dann
  \begin{equation*}
    \widehat{X \otr M} \longrightarrow \widehat{X \otr P}
  \end{equation*}
  ein zerfallender Monomorphismus.
\end{Folgerung}

\begin{Beweis}
  Wie in Beispiel~1 ist $\widehat{M} \to \widehat{P}$
  ein zerfallender Monomorphismus, also auch $X \otr \widehat{M} \to X \otr
  \widehat{P}$. Weil auch $\widehat{M}$ und $\widehat{P}$ totalsepariert
  sind, liefert (\ref{1.1}) die Behauptung.
\end{Beweis}

\begin{Folgerung}\label{1.9}
  Sei $f_i \colon F_i \twoheadrightarrow M_i$ eine Familie von Epimorphismen
  $(i \in I)$, bei der jedes $M_i$ separiert und jedes $F_i$ flach ist.
  Außerdem gelte für jedes $n \geq 1$, dass fast alle $\overline{f_i} \colon
  F_i / \mathfrak{m}^{n} F_i \to M_i / \mathfrak{m}^{n} M_i$ Isomorphismen
  sind. Dann ist
  \begin{equation*}
    \widehat{\coprod M_i} \longrightarrow \widehat{\prod M_i}
  \end{equation*}
  ein zerfallender Monomorphismus.
\end{Folgerung}

\begin{Beweis}
  Mit $M =  \coprod M_i$ und $P = \prod M_i$ ist nach (\ref{1.6}) nur noch
  zu zeigen, dass alle $P/M + \mathfrak{m}^{n}P$ als
  $R/\mathfrak{m}^{n}$-Moduln frei sind. Mit $F = \coprod F_i$ und $Q =
  \prod F_i$ ist auch $Q/F$ flach, also $Q/F + \mathfrak{m}^{n}Q$ als
  $R/\mathfrak{m}^{n}$-Modul frei. Bleibt zu zeigen, dass der Epimorphismus
  $g = \prod f_i \colon Q \to P$ einen Isomorphismus $g_n \colon Q/F +
  \mathfrak{m}^{n}Q \to P/M + \mathfrak{m}^{n}P$ induziert: Aus $u = (u_i)
  \in Q$ und $g(u) \in M + \mathfrak{m}^{n}P$, d.\,h. $f_i(u_i) \in
  \mathfrak{m}^{n}M_i$ für fast alle $i$, folgt nach Voraussetzung $u_i \in
  \mathfrak{m}^{n}F_i$ für fast alle $i$, d.\,h. $u \in F + \mathfrak{m}^{n}Q$.
\end{Beweis}

\begin{Satz}\label{1.10}
  Seien in $M = \coprod M_i$ alle $M_i$ endlich erzeugt und gelte für jedes
  $n \geq 1$, dass fast alle $M_i/\mathfrak{m}^{n} M_i$ als
  $R/\mathfrak{m}^{n}$-Moduln frei sind. Dann ist $\widehat{M}$ rein-injektiv.
\end{Satz}

\begin{Beweis}
  Weil alle $\widehat{M_i}$ rein-injektiv sind, ist es auch $\prod
  \widehat{M_i} \cong \widehat{\prod M_i}$, und wir müssen nur noch die
  Voraussetzungen in (\ref{1.9}) nachweisen. Jedes $M_i$ besitzt, weil
  endlich erzeugt, eine projektive Hülle $f_i \colon F_i \twoheadrightarrow
  M_i$, und dann ist auch $\overline{f_i}\colon F_i/\mathfrak{m}^{n}F_i \to
  M_i/\mathfrak{m}^{n}M_i$ eine projektive Hülle über $R/\mathfrak{m}^{n}$.
  War nun $M_i/\mathfrak{m}^{n}M_i$ über $R/\mathfrak{m}^{n}$ frei, ist
  $\overline{f_i}$ ein Isomorphismus.
\end{Beweis}

\begin{Folgerung}\label{1.11}
  Sei $M = \coprod_{i=1}^{\infty} R/\mathfrak{a}_i$ und gelte für jedes $n
  \geq 1$, dass $\mathfrak{m}^n$ fast alle $\mathfrak{a}_i$ enthält. Dann
  ist $\widehat{M}$ rein-injektiv.
\end{Folgerung}

\begin{Beweis}
  Mit $M_i = R/\mathfrak{a}_i$ folgt aus $\mathfrak{m}^n \supset
  \mathfrak{a}_i$, dass $M_i / \mathfrak{m}^n M_i \cong R / \mathfrak{m}^n$ ist.
\end{Beweis}

Die Bedingung in (\ref{1.11}) ist z.\,B. dann erfüllt, wenn $\mathfrak{a}_i
= \mathfrak{m}^i$ ist für alle $i \geq 1$, oder wenn $R$ vollständig und
$\mathfrak{a}_1 \supset \mathfrak{a}_2 \supset \mathfrak{a}_3 \supset
\ldots$ eine absteigende Folge von Idealen ist mit $\bigcap_{i=1}^{\infty}
\mathfrak{a}_i = 0$ (Theorem von Chevalley).

\section{Wann ist $\widehat{M}$ totalsepariert?}

Ist $M$ totalsepariert, so ist die in (\ref{1.1}) angegebene Isomorphie
sogar äquivalent damit, dass $\widehat{M}$ totalsepariert ist. Wir wissen
nicht, ob das immer der Fall ist, geben aber in (\ref{2.3}) eine Reihe von
dazu äquivalenten Bedingungen an. Insbesondere untersuchen wir den
Zusammenhang zwischen totalseparierten Moduln und strikt rein-wesentlichen
Erweiterungen. Dabei heißt eine Modulerweiterung $A \subset B$ \emph{strikt
  rein-wesentlich} (\cite{012} p.~8), wenn $A \subset B$ rein-wesentlich ist
und aus $B \subset Y$, $A$ rein in $Y$ stets folgt $B$ rein in $Y$. Z.\,B.
gilt für jeden Basis-Untermodul $F$ eines flachen $R$-Moduls $M$: Genau dann
ist $F$ rein-wesentlich (strikt rein-wesentlich) in $M$, wenn $M$ separiert
(totalsepariert) ist.

\bigskip

Die im dritten Abschnitt von \cite{012} aufgestellten Rechenregeln für
strikt rein-wesentliche Erweiterungen ergänzen wir durch

\begin{Lemma}\label{2.1}
  Sei $A \subset B$ strikt rein-wesentlich und $A \subset B_1 \subset B$
  ein Zwischenmodul, so dass $B_1$ rein in $B$ ist. Dann sind auch $A
  \subset B_1$ und $B_1 \subset B$ strikt rein-wesentlich.
\end{Lemma}

\begin{Beweis}
  (a) $A \subset B_1$ ist rein-wesentlich, denn aus $X \subset B_1$, $X \cap
  A = 0$ und $(X \oplus A) / X$ rein in $B_1 / X$ folgt $(X \oplus A) / X$
  rein in $B/X$, also $X=0$.
  (b) $A \subset B_1$ ist strikt rein-wesentlich, denn aus $B_1 \subset Y$
  und $A$ rein in $Y$ folgt, dass in der Fasersumme
  \begin{equation*}
    \xymatrix{
      B \ar@{}[d]|{\textstyle\cup} \ar@{-->}[r]^h & Y' \\
      B_1\ar@{}[r]|{\textstyle\subset} & Y \ar@{-->}[u]_g
    }    
  \end{equation*}
  auch $g$ ein reiner Monomorphismus ist, also auch $A \subset B
  \xrightarrow{h}Y'$. Weil $A \subset B$ strikt rein-wesentlich war, folgt
  die Reinheit von $h$, also auch von $B_1 \subset Y$.
  (c) $B_1 \subset B$ ist rein-wesentlich, denn aus $X \subset B$, $X \cap
  B_1 = 0$ und $(X \oplus B_1) / X$ rein in $B/X$ folgt $X \cap A = 0$ und
  $(X \oplus A) / X$ rein in $B/X$, also $X=0$. (d) $B_1 \subset B$ ist
  strikt rein-wesentlich, denn aus $B \subset Y$ und $B_1$ rein in $Y$ folgt
  $A$ rein in $Y$, also nach Voraussetzung $B$ rein in $Y$.
\end{Beweis}

\begin{Folgerung}\label{2.2}
  Sei $M \subset Y$ eine beliebige Modulerweiterung.
  \begin{abc}
  \item Ist $M$ totalsepariert und $M \subset Y$ strikt rein-wesentlich, so
    ist auch $Y$ totalsepariert.
  \item Ist $Y$ totalsepariert und $M \subset Y$ rein und dicht, so ist $M
    \subset Y$ strikt rein-wesentlich.
  \end{abc}
\end{Folgerung}

\begin{Beweis}
  (a) Die rein-injektive Hülle $N$ von $M$ ist ebenfalls totalsepariert, und
  weil man $Y \subset N$ wählen kann, ist auch $Y$ rein in $N$ und deshalb
  totalsepariert.\\
  (b) Nach (\cite{011} Lemma~1.9) ist $M \subset Y$ rein-wesentlich, so dass
  man wieder $Y \subset N$ wählen kann, und im Diagramm
  \begin{equation*}
    \xymatrix{
      \prod_{n\geq 1} M/\mathfrak{m}^nM \ar[r]^{\cong} & \prod_{n\geq 1}
      Y/\mathfrak{m}^nY \ar[r]^{\beta}
      & \prod_{n\geq 1} N/\mathfrak{m}^nN \\
      & Y\ar[u]^{\alpha}\ar@{}[r]|{\textstyle\subset} & N\ar[u]
    }    
  \end{equation*}
  sind dann $\alpha$ und $\beta$ reine Monomorphismen. Damit ist auch $Y$
  rein in $N$ und deshalb $M \subset Y$ nach (\ref{2.1}) strikt rein-wesentlich.
\end{Beweis}

\begin{Satz}\label{2.3}
  Für einen totalseparierten $R$-Modul $M$ sind äquivalent:
  \begin{iii}
    \item $\widehat{M}$ ist totalsepariert.
    \item $M \subset \widehat{M}$ ist strikt rein-wesentlich.
    \item Für jeden endlich erzeugten $R$-Modul $X$ ist
      \begin{equation*}
        \widehat{X \otr M} \cong X \otr \widehat{M}.
      \end{equation*}
    \item Ist $\alpha\colon M \to Y$ ein reiner Monomorphismus und $Y$
      separiert, so ist auch\\ $\widehat{\alpha}\colon \widehat{M} \to
      \widehat{Y}$ ein reiner Monomorphismus.
  \end{iii}
\end{Satz}

\begin{Beweis}
  (i $\leftrightarrow$ ii) Nach Voraussetzung ist $M$ rein in $\widehat{M}$
  (und natürlich immer dicht), so dass (\ref{2.2}) mit $Y = \widehat{M}$ die
  Äquivalenz liefert.\\
  Bei (i $\to$ iii) ist nach (\ref{1.1}) sogar die kanonische Abbildung
  $\widehat{X \otr M} \to X \otr \widehat{M}$ ein Isomorphismus, und bei
  (iii $\to$ i) ist $X \otr \widehat{M}$ separiert für jeden endlich
  erzeugten $R$-Modul $X$, d.\,h. $\widehat{M}$ totalsepariert.\\
  Bei (i $\to$ iv) ist in $\widehat{M} \subset \prod_{n\geq 1}
  M/\mathfrak{m}^nM \xrightarrow{\cong} \prod_{n\geq 1}
  \widehat{M}/\mathfrak{m}^n \widehat{M}$ nach Voraussetzung die erste
  Inklusion rein, im kommutativen Diagramm
  \begin{equation*}
    \xymatrix@C=4em@R=4em{
      \prod_{n\geq 1} M/\mathfrak{m}^nM \ar[r]^{\prod
        \alpha_n}\ar@{}[d]|{\textstyle \cup} & \prod_{n\geq 1}
      Y/\mathfrak{m}^nY\ar@{}[d]|{\textstyle \cup}
    \\
    \qquad\widehat{M}\qquad\ar[r]^{\widehat{\alpha}} & \qquad\widehat{Y}\qquad
    }
  \end{equation*}
  aber auch die obere Zeile, also auch $\widehat{\alpha}$. Bei (iv $\to$ i)
  gilt speziell für die rein-injektive Hülle $N$ von $M$, dass $N \cong
  \widehat{N}$ totalsepariert ist, also auch der nach Voraussetzung reine
  Untermodul $\widehat{M}$.
\end{Beweis}

\begin{Bemerkung}\label{2.4}
  Ist $M$ nur separiert, aber rein in $\widehat{M}$, kann man die
  Eigenschaften "`$\widehat{M}$ ist totalsepariert"' bzw. "`$\widehat{M}$
  ist rein-injektiv"' auch an der rein-injektiven Hülle $N$ von $M$ ablesen:
  Ist $N_1/M$ der größte radikalvolle Untermodul von $N/M$, gilt nach
  (\cite{011} Folgerung~1.6) $H(N_1) = H(N)$ sowie $\widehat{M} \cong N_1 /
  H(N_1)$ (mit $H(X) = \bigcap_{n \geq 1} \mathfrak{m}^n X$ für jeden
  $R$-Modul $X$), und damit folgt:
  \begin{center}
    $\widehat{M}$ totalsepariert $\iff$ $N_1/M$ rein in $N/M$;\\[1ex]
    $\widehat{M}$ rein-injektiv $\iff$ $N/M$ ist radikalvoll.
  \end{center}
  Mit Hilfe des nächsten Lemmas wollen wir eine Reihe von nicht
  totalseparierten Moduln und nicht strikt rein-wesentlichen Erweiterungen
  konstruieren:
\end{Bemerkung}

\begin{Lemma}\label{2.5}
  Sei $M \subset A$ eine Modulerweiterung, so dass $A$ separiert, $M$ rein
  und dicht in $A$ ist. Sei $\mathfrak{c}$ ein Ideal in $R$ mit der
  Eigenschaft (*) $A/M$ ist nicht $\mathfrak{c}$-teilbar und $\mathfrak{c}
  \not\subset \Ass(A)$.\\
  Dann ist $B := M + \mathfrak{c} \cdot A$ ein Modul zwischen $M$ und $A$,
  so dass $M$ auch rein und dicht in $B$ ist, aber $B$ nicht rein in $A$.
  Insbesondere ist $M$ nicht strikt rein-wesentlich in $B$ und deshalb $B$
  nicht totalsepariert.
\end{Lemma}

\begin{Beweis}
  Mit $A/M$ ist auch $\mathfrak{c} \cdot A/M = B/M$ radikalvoll, also $M$
  rein und dicht in $B$ und deshalb nach (\cite{011} Lemma~1.9)
  rein-wesentlich. Wegen (*) gibt es ein $r \in \mathfrak{c}$, das NNT bzgl.
  $A$ ist, und aus $B \subsetneqq A$ folgt $rB \subsetneqq rA$, $rB
  \subsetneqq B \cap rA$, so dass $B$ nicht rein in $A$ ist. Damit ist der
  Zusatz klar. $B$ ist nicht einmal totalreduziert (siehe \cite{011}
  Folgerung 1.10), denn $B/rB$ hat den radikalvollen Untermodul $rA/rB \neq 0$.
\end{Beweis}

\begin{Bemerkung}\label{2.6}
  Das Ideal $\mathfrak{c}$ mit (*) kann man variieren: Vergrößern bis zu
  einen $\mathfrak{q} \in \operatorname{Koass}(A/M)$ oder verkleinern bis zu
  einem 
  $(s) \subset R$ mit $\operatorname{Ann}_A(s) = 0$. Genau dann gibt es also
  ein Ideal $\mathfrak{c}$ mit (*), wenn es ein Ringelement $s$ gibt, das
  auf $A$ injektiv und auf $A/M$ nicht surjektiv operiert, d.\,h. wenn
  $\bigcup \operatorname{Koass}(A/M) \not\subset \bigcup \Ass(A)$ ist.
\end{Bemerkung}

\noindent
\textbf{Beispiel 1} Sei $A$ separiert und flach $\neq 0$, $M$ ein
Basis-Untermodul von $A$. Genau dann gibt es ein Ideal $\mathfrak{c}$ mit
(*), wenn $A/M$ \emph{nicht} teilbar ist. (Das ist z.\,B. der Fall, wenn $R$
ein Integritätsring der Dimension $\geq 2$ ist und $M \cong R^{(\mathbb{N})}$, $A =
\widehat{M}$.)

\bigskip

\noindent
\textbf{Beispiel 2} Sei $M$ endlich erzeugt und $A = \widehat{M}$. Genau
dann gibt es ein Ideal $\mathfrak{c}$ mit (*), wenn es einen NNT $s$ bzgl.
$M$ gibt, so dass $R/\operatorname{Ann}_R(M) + (s)$ \emph{nicht} vollständig
ist. (Das ist z.\,B. der Fall, wenn $R$ ein abzählbarer Integritätsring der
Dimension $\geq 2$ ist und $M \cong R^n$, $n \geq 1$.)

\bigskip

\noindent
\textbf{Beispiel 3} Sei $M = X^{(I)}$ mit $X$ endlich erzeugt $\neq 0$ und
$I$ unendlich, $A = \widehat{M}$. Genau dann gibt es ein Ideal
$\mathfrak{c}$ mit (*), wenn es einen NNT $s$ bzgl. $X$ gibt, so dass
$R/\operatorname{Ann}_R(M) +(s)$ \emph{nicht} artinsch ist. (Denn nach
(3.8') ist $\operatorname{Koass}(A/M) = \{\mathfrak{p} \in
\operatorname{Spec}(R) \mid \operatorname{Ann}_R(X) \subset \mathfrak{p}
\subsetneqq \mathfrak{m}\}$.

\section{Wann ist $A$ rein in $P$?}

Stets sei in diesem Abschnitt $(M_i \mid i \in I)$ eine Familie von
$R$-Moduln, $M = \coprod M_i$, $P = \prod M_i$ und $A = \bigcap_{n \geq 1}(M
+ \mathfrak{m}^n P)$ der $\mathfrak{m}$-adische Abschluss von $M$ in $P$.
Auch jetzt fragen wir uns, wann $A$ totalsepariert ist. Nach (\ref{2.2}) ist
das genau dann der Fall, wenn $M$ totalsepariert ist und $A$ rein in $P$.
Die zweite Bedingung wollen wir genauer untersuchen.

\begin{Satz}\label{3.1}
  Sei $M = \coprod M_i$ beliebig und existiere zum Ideal $\mathfrak{b}$ in
  $R$ ein $e \geq 1$ mit $\mathfrak{m}^f M\, \cap\, \mathfrak{b} M \subset
  \mathfrak{m}^{f-e}\mathfrak{b} M$ für alle $f \geq e$. Dann gilt
  \begin{equation*}
    A \cap \mathfrak{b} P \subset \mathfrak{b} A + H(P)\ .
  \end{equation*}
\end{Satz}

\begin{Beweis}
  Sei im \textbf{1. Schritt} $M$ separiert, d.\,h. $H(M) = \bigcap
  \mathfrak{m}^n M = 0$. Dann müssen wir $A \cap \mathfrak{b} P =
  \mathfrak{b} A$ zeigen und gehen wie in (\cite{008} Proposition 2.1.9)
  vor: Zu $(x_i) \in A \cap \mathfrak{b} P$ gibt es nach Definition eine
  endliche Teilmenge $I_0$ von $I$ mit $x_i \in \mathfrak{m}^e M_i$ für alle
  $i \notin I_0$. Sei $\mathfrak{b} = (b_1,\dotsc,b_m)$.\\
  \emph{Beh. 1} Für alle $i \notin I_0$ gibt es eine Darstellung $x_i =
  \sum_{\mu = 1}^m b_{\mu} u_{i\mu}$ mit $u_{i1}, \dotsc, u_{im}\in M_i$ und
  der Eigenschaft
  \begin{equation*}
    (\Delta)\quad n \geq 1 \quad\text{und}\quad x_i \in
    \mathfrak{m}^{n+e}M_i \implies u_{i1}, \dotsc, u_{im}\in \mathfrak{m}^n
    M_i\ .
  \end{equation*}
  Klar ist $x_i = 0$: Setze $u_{i1} = \cdots = u_{im} = 0$. Bei $x_i \neq 0$
  gibt es, weil $M_i$ separiert ist, genau ein $t_i \geq e$ mit $x_i \in
  \mathfrak{m}^{t_i} M_i \setminus \mathfrak{m}^{t_i+1} M_i$, aus $x_i \in
  \mathfrak{m}^{t_i} M_i \cap \mathfrak{b}M_i$ folgt nach Voraussetzung $x_i
  \in \mathfrak{m}^{t_i-e}\mathfrak{b}M_i$, d.\,h. $x_i = \sum_{\mu = 1}^m
  b_{\mu} u_{i\mu}$ mit $u_{i\mu} \in \mathfrak{m}^{t_i-e} M_i$. Die
  Eigenschaft $(\Delta)$ folgt dann aus der Maximalität von $t_i$.\\
  \emph{Beh. 2} Für alle $i \in I_0$ setze $y_i = u_{i1} = \dotsc = u_{im} =
  0$. Für alle $i \notin I_0$ setze $y_i = x_i$, und dann ist $(y_i) =
  \sum_{\mu=1}^m b_{\mu} (u_{i\mu})$, $(y_i)- (x_i) \in M \cap \mathfrak{b}P
  = \mathfrak{b}M$, so dass $(y_i) \in \mathfrak{b}A$, d.\,h. $(u_{i1}),
  \dotsc, (u_{im}) \in  A$ zu zeigen bleibt: Zu festem $n \geq 1$ gibt es
  nach Definition von $A$ eine endliche Teilmenge $I_1 \supset I_0$ von $I$
  mit $x_i \in \mathfrak{m}^{n+e} M_i$ für alle $i \notin I_1$, mit
  $(\Delta)$ folgt dann $u_{i1},\dotsc,u_{im} \in \mathfrak{m}^n M_i$ für
  alle $i \notin I_1$ wie verlangt.\\[1ex]
  Sei im \textbf{2. Schritt} $M$ beliebig und $\overline{P} = P/H(P)$. Im
  kommutativen Diagramm
  \begin{equation*}
    \xymatrix{
      \overline{M}\ar[d]\ar@{}[r]|{\textstyle \subset}  &
      \overline{A}\ar[d]\ar@{}[r]|{\textstyle \subset} & 
      \overline{P}\ar[d]_{\omega}^{\cong} \\
      \coprod \overline{M}_i\ar@{}[r]|{\textstyle \subset} &
      \omega(\overline{A})\ar@{}[r]|{\textstyle \subset} & \prod
      \overline{M}_i 
    }
  \end{equation*}
  ist dann $\overline{A} = A/H(P)$ der $\mathfrak{m}$-adische Abschluss von
  $\overline{M}$ in $\overline{P}$ und, wie eine elementare Rechnung zeigt,
  auch $\mathfrak{m}^f \overline{M} \cap \mathfrak{b}\overline{M} \subset
  \mathfrak{m}^{f-e}\mathfrak{b} \overline{M}$ für alle $f \geq e$. Nach dem
  ersten Schritt folgt $\omega(\overline{A}) \cap \mathfrak{b}(\prod
  \overline{M}_i) = \mathfrak{b} \cdot \omega(\overline{A})$, d.\,h.
  $\overline{A} \cap \mathfrak{b}\overline{P} = \mathfrak{b}\overline{A}$,
  d.\,h. $(A \cap \mathfrak{b}P) + H(P) = \mathfrak{b}A + H(P)$ wie
  behauptet.
\end{Beweis}

\begin{Folgerung}\label{3.2}
  Sei $M = \coprod M_i$ beliebig und $\mathfrak{b} \subset R$ ein
  \emph{offenes} Ideal. Dann gilt
  \begin{equation*}
    A \cap \mathfrak{b}P = \mathfrak{b}A + H(P) \quad\text{und}\quad A = M +
    \mathfrak{b}A + H(P)\ .
  \end{equation*}
\end{Folgerung}

\begin{Beweis}
  Aus $\mathfrak{m}^e \subset \mathfrak{b}$ für ein $e \geq 1$ folgt
  $\mathfrak{m}^fM \subset \mathfrak{m}^{f-e}\mathfrak{b}M$ für alle $f \geq
  e$. Also ist die Voraussetzung in (\ref{3.1}) erfüllt, wegen $H(P) \subset
  \mathfrak{b}P$ gilt im Satz $A \cap \mathfrak{b}P = \mathfrak{b}A + H(P)$,
  und Addition mit $M$ liefert die zweite Gleichung.
\end{Beweis}

Aus ihr folgt sofort: Ist $M = \coprod M_i$ separiert, also $H(M) = 0$,
$H(P) = 0$, $A = M + \mathfrak{m}A$, so ist $A/M$ radikalvoll. Aber die
Bedingung "`$M$ separiert"' ist nicht notwendig.

\begin{Folgerung}
  Ist $M = \coprod_{i=1}^{\infty} M_i$, so ist $A/M$
  radikalvoll. 
\end{Folgerung}

\begin{Beweis}
  Nach (\ref{3.2}) genügt es, $H(P) \subset \mathfrak{m}A$ zu zeigen: $x =
  (x_i) \in H(P) = \prod H(M_i) \Rightarrow$ mit $\mathfrak{m}= (b_1,
  \dotsc, b_m)$ gilt $x_i \in  \mathfrak{m}^i M_i$, $x_i =\sum_{\mu=1}^m
  b_{\mu} u_{i\mu}$ mit $u_{i1},\dotsc, u_{im} \in \mathfrak{m}^{i-1}M_i
  \Rightarrow (x_i) = \sum_{\mu=1}^m b_{\mu} (u_{i\mu})$, alle $(u_{i\mu})
  \in A$, also $x \in \mathfrak{m}A$.
\end{Beweis}

\begin{Folgerung}\label{3.3}
  Für $M = \coprod M_i$ sind äquivalent:
  \begin{iii}
    \item $A/M$ ist radikalvoll.
    \item $A \cap \mathfrak{m}P = \mathfrak{m}A$.
    \item $A \cap \mathfrak{b}P = \mathfrak{b}A$ für jedes offene Ideal
      $\mathfrak{b}$.
    \item $H(A) = H(P)$
  \end{iii}
\end{Folgerung}

\begin{Beweis}
  (i $\to$ iii) Aus $M + \mathfrak{m}A = A$ und $\mathfrak{m}^e \subset
  \mathfrak{b}$ für ein $e \geq 1$ folgt $M + \mathfrak{b}A = A$,
  $\mathfrak{b}A + (M \cap \mathfrak{b}P) = A \cap \mathfrak{b}P$, und
  natürlich ist $M \cap \mathfrak{b}P = \mathfrak{b}M$ enthalten in
  $\mathfrak{b}A$. (iii $\to$ ii) ist klar, und bei (ii $\to$ i) ist $M +
  \mathfrak{m}A = (M + \mathfrak{m}P) \cap A = A$.\\
  Ebenso ist (iii $\to$ iv) klar, und bei (iv $\to$ i) genügt $H(P) \subset M
  + \mathfrak{m}A$: Nach (\ref{3.2}) folgt dann $A = M + \mathfrak{m}A$.
\end{Beweis}

\begin{Folgerung}\label{3.5neu}
  Sind in $M = \coprod M_i$ fast alle $H(M_i)$ radikalvoll,
  so ist auch $A/M$ radikalvoll.
\end{Folgerung}

\begin{Beweis}
  Aus der Voraussetzung folgt, dass $\prod H(M_i) / \coprod H(M_i) = H(P) /
  H(M) \cong (M + H(P))/M$ radikalvoll ist, also $H(P) \subset M +
  \mathfrak{m}A$, also wieder mit (\ref{3.2}) die Behauptung.
\end{Beweis}

Sei $M = \coprod M_i$ wieder beliebig und $(x_i) \in A$. Für jedes $n \geq
1$ ist dann die Menge $I_n = \{i \in I \mid x_i \notin \mathfrak{m}^n M_i \}$
endlich und deshalb $\bigcup_{n\geq 1} I_n = \{i \in I \mid x_i \notin
H(M_i) \}$ abzählbar. War also $I$ überabzählbar, gibt es ein $i_0 \in I$
mit $x_{i_0} \in H(M_{i_0})$. Das verwenden wir im folgenden

\bigskip

\noindent
\textbf{Beispiel}\hspace{0.2em} Sei $R$ ein diskreter Bewertungsring und $M =
X^{(I)}$ mit $H(X)$ nicht radikalvoll, $I$ überabzählbar. Dann ist $A/M$
\emph{nicht} radikalvoll.

\begin{Beweis}
  Sei $\mathfrak{m} = (p)$, $u \in H(X) \setminus p H(X)$ und $x = (x_i)$
  mit $x_i = u$ für alle $i$. Klar ist $x \in H(P)$, und wäre $A/M$
  radikalvoll, also nach (\ref{3.3}) $H(P) \subset pA$, folgte $x = py$ für
  ein $y = (y_i) \in A$. Nach der Vorbemerkung gilt $y_{i_0} \in H(X)$ für
  ein $i_0 \in I$, also $u = p y_{i_0} \in p H(X)$ entgegen der Wahl.
  (Bekanntlich gibt es zu jedem $R$-Modul $U$ eine Erweiterung $U \subset X$
  mit $H(X) = U$.)
\end{Beweis}

\begin{Satz}\label{3.4}
  Sei $M = \coprod M_i$ und $M/H(M)$ flach. Dann ist auch $P/A$ flach.
\end{Satz}

\begin{Beweis}
  Falls alle $M_i = R$ sind, ist das wohlbekannt (Raynaud-Gruson \cite{004}
  p.~77). Aber auch in unserer Situation folgt mit $\overline{P} = P/H(P)$
  und dem Diagramm aus dem zweiten Beweisschritt von (\ref{3.1}), dass man
  $H(P) = 0$, d.\,h. $H(M) = 0$ annehmen kann. Nach Artin-Rees gibt es zu
  \emph{jedem} Ideal $\mathfrak{b}$ ein $e \geq 1$ mit $\mathfrak{m}^f \cap
  \mathfrak{b} \subset \mathfrak{m}^{f-e} \mathfrak{b}$ für alle $f \geq
  e$, weil $M$ flach ist, folgt $\mathfrak{m}^fM  \cap \mathfrak{b}M \subset
  \mathfrak{m}^{f-e} \mathfrak{b} M$ für dieselben $f$, also nach
  (\ref{3.1}) $A \cap \mathfrak{b}P = \mathfrak{b}A$. Weil auch $P$ flach
  ist, folgt $\operatorname{Tor}_1^R (R/\mathfrak{b}, P/A) = 0$ für jedes
  Ideal $\mathfrak{b}$, also die Behauptung.
\end{Beweis}

\bigskip
\noindent
\textbf{Beispiele} (1) Sei $M = X^{(I)}$ mit $X$ endlich erzeugt. Dann ist
$A$ rein in $P$, denn mit den Bezeichnungen wie in (\cite{004} p.~77) ist $A
= C(I,X) \cong C(I,R) \otr X$, und weil $C(I,R)$ rein in $R^I$ ist, ist es
auch $A$ in $P$.\\[1ex]
(2) Sei $R$ ein 1-dim. Integritätsring und $M = \coprod M_i$ torsionsfrei.
Dann ist $A$ rein in $P$, denn alle $H(M_i)$ sind teilbar, also radikalvoll,
so dass nach (\ref{3.5neu}) $A/M$ radikalvoll, ja sogar teilbar und
torsionsfrei ist und aus $A/M \subset^{\oplus} P/M$ die Behauptung
folgt.\\[1ex]
(3) Sei $f_i \colon F_i \twoheadrightarrow M_i$ eine Familie von
Epimorphismen wie in (\ref{1.9}), in der auch jedes $F_i$ separiert ist. Mit
$B = \bigcap_{n\geq 1} (F + \mathfrak{m}^n Q)$ ist dann nach (\ref{3.4})
auch $Q/B$ flach, außerdem im Diagramm
\begin{equation*}
  \xymatrix@C=4em@R=4em{
    Q/B\ \ar@{^{(}->}[r] \ar@{->>}[d]_{\overline{g}} & \prod_{n \geq 1}Q/F + \mathfrak{m}^n Q \ar[d]^{\prod g_n}\\
    P/A\ \ar@{^{(}->}[r] & \prod_{n \geq 1}P/M + \mathfrak{m}^n P
  }
\end{equation*}
nach dem Beweis von (\ref{1.9}) jedes $g_n$ ein Isomorphismus, also auch
$\overline{g}$, und damit ist $P/A$ flach.\hfill$\Box$

\bigskip

Unter den Voraussetzungen von (\ref{3.4}) ist $A/M$ radikalvoll, d.\,h.
$\mathfrak{m} \notin \operatorname{Koass}(A/M)$. Ist aber $M$ groß genug,
sind alle Primideale $\neq \mathfrak{m}$ zu $A/M$ koassoziiert:

\begin{Satz}\label{3.5}
  Sei $M =  \coprod M_i$ und $M/H(M)$ flach, und seien unendlich viele $M_i$
  \emph{nicht} radikalvoll. Dann ist
  \begin{equation*}
    \operatorname{Koass}(A/M) = \operatorname{Spec}(R) \setminus
    \{\mathfrak{m}\}\,.
  \end{equation*}
\end{Satz}

\begin{Beweis}
  Mit $\overline{P} = P/H(P)$ hat man den Epimorphismus $A/M
  \twoheadrightarrow \overline{A}/\overline{M}$, kann also mit Hilfe des
  Diagramms in (\ref{3.1}) gleich $H(P) = 0$ annehmen, d.\,h. $M = \coprod
  M_i$ separiert und flach, unendlich viele $M_i \neq 0$. Weil dann $A/M$
  radikalvoll und flach ist, gilt nach (\cite{010} Folgerung~2.3)
  $\operatorname{Koass}(A/M) = \{\mathfrak{p} \in \operatorname{Spec}(R)
  \mid \mathfrak{p} \cdot A/M \neq A/M \}$ und wir müssen deshalb für jedes
  Primideal $\mathfrak{q} \subsetneqq \mathfrak{m}$ zeigen, dass $B = M +
  \mathfrak{q}A$ echt kleiner als $A$ ist:

  \bigskip
  \noindent
  Sind die $M_{j_0}, M_{j_1}, M_{j_2}, \dotsc$ ungleich null, definieren wir
  mit $I'=\{j_0,j_1,j_2,\dotsc \}$ (paarw. versch.) und $I = I' \dotcup I''$
  ein Element $x = (x_i) \in P$ durch $x_i \in \mathfrak{m}^m M_i \setminus
  (\mathfrak{q} \cap \mathfrak{m}^m)M_i$ falls $i \in I'$, $i =j_m$ (das
  geht, weil $\mathfrak{q} \cap \mathfrak{m}^m \subsetneqq \mathfrak{m}^m$
  ist und $M_i$ einen reinen Untermodul $\cong R$ besitzt), $x_i = 0$ falls
  $i \in I''$. Es ist dann $x \in A$, denn für festes $n \geq 1$ ist $\{i
  \in I \mid x_i \notin  \mathfrak{m}^nM_i\} \subset
  \{j_0,j_1,\dotsc,j_{n-1} \}$ endlich. Wäre aber $x \in B$, d.\,h. $x-y \in
  \mathfrak{q}A$ mit $y \in M$, folgte $\{i \in I \mid x_i \notin
  \mathfrak{q} M_i \} \subset I'$, worin die erste Menge endlich und die
  zweite unendlich ist, also $x_i \in \mathfrak{q}M_i$ für ein $i = j_m$.
  Weil $M_i$ flach war, folgt $x_i \in (\mathfrak{q} \cap
  \mathfrak{m}^m)M_i$ und das ist unmöglich.
\end{Beweis}

\begin{Folgerung}\label{3.6}
  Sei $M = X^{(I)}$ mit $X$ endlich erzeugt $\neq 0$ und $I$ unendlich. Dann ist
  \begin{equation*}
    \operatorname{Koass}(A/M) = \{\mathfrak{p} \in \operatorname{Spec}(R)
    \mid \operatorname{Ann}_R(X) \subset \mathfrak{p} \subsetneqq
    \mathfrak{m} \ \}\,.
  \end{equation*}
\end{Folgerung}

\begin{Beweis}
  Mit $A = C(I,X)$ wie in (\cite{004} p.~77) folgt aus
  \begin{equation*}
    \xymatrix@C=3em@R=3em{
      R^{(I)} \otr X \ar[r]\ar[d]^{\cong} & C(I,R) \otr X\ar[d]^{\cong} \\
      M \ar@{}[r]|{\textstyle \subset} & A
    }
  \end{equation*}
  und $Y := C(I,R)/R^{(I)}$, dass $A/M \cong Y \otr X$ ist, mit der
  entsprechenden Formel für $\operatorname{Koass}$ in (\cite{009}
  Folgerung~3.2) also $\operatorname{Koass}(Y) \cap \operatorname{Supp}(X) =
  \{\mathfrak{p} \in \operatorname{Spec}(R) \mid \mathfrak{p} \neq
  \mathfrak{m} \text{ und } \operatorname{Ann}_R(X) \subset \mathfrak{p}
  \}$.
\end{Beweis}

Wir wollen die beiden letzten Aussagen noch für $\widehat{M}$ statt für $A$
fomulieren. Ist $M = \coprod M_i$ separiert, induzieren die kanonischen
Abbildungen $\beta_n = A \subset M + \mathfrak{m}^n P \to (M +
\mathfrak{m}^n P)/\mathfrak{m}^n P \xrightarrow{\cong} M/\mathfrak{m}^n M$
einen Monomorphismus $\beta\colon A \to \widehat{M}$ mit $\operatorname{Bi}
\beta + \mathfrak{m}^n \widehat{M} = \widehat{M}$, und in
\begin{equation*}
  \xymatrix@C=3.5em@R=3.5em{
    M \ar@{}[r]|{\textstyle\subset}\ar[d]_{\mu} &
    A\ar[dl]_{\beta}\ar[d]^{\alpha} \\
    \widehat{M}\ar[r]^{\cong} & \widehat{A}
  }
\end{equation*}
sind beide Teildreiecke kommutativ. Ist $A$ sogar totalsepariert, wird
$\alpha$, also auch $\beta$ ein reiner Monomorphismus, ebenso die induzierte
Abbildung $\overline{\beta}\colon A/M \to \widehat{M}/
\operatorname{Bi}\mu$.

\bigskip

\noindent
\textbf{Satz 3.7'}
  {\itshape Sei $M = \coprod M_i$ totalsepariert und flach, und seien unendlich viele
  $M_i \neq 0$. Dann ist $\operatorname{Koass}(\widehat{M}/M) =
  \operatorname{Spec}(R) \setminus \{\mathfrak{m}\}$.}

\bigskip
\noindent
\textbf{Folgerung 3.8'}
  {\itshape Sei $M = X^{(I)}$ mit $X$ endlich erzeugt $\neq 0$ und $I$ unendlich. Dann
  ist $\operatorname{Koass}(\widehat{M}/M) =  \{\mathfrak{p} \in
  \operatorname{Spec}(R) \mid \operatorname{Ann}_R(X) \subset \mathfrak{p}
  \subsetneqq \mathfrak{m} \}$.}

\bigskip

\begin{Beweis}
  Beide Male ist $A$ rein in $P$, also auch $A$ totalsepariert und deshalb
  $\overline{\beta}\colon A/M \to \widehat{M}/M$ ein reiner Monomorphismus.
  Es folgt $\operatorname{Koass}(\widehat{M}/M) = \operatorname{Koass}(A/M)$
  wie gewünscht.
\end{Beweis}

\end{document}